\documentclass[12pt]{amsart}

\newcommand\I{\mathrm I}
\newcommand\R{\mathrm R}
\newcommand\nn{\mathrm n}

\renewcommand\c{\mathfrak c}
\newcommand\g{\mathfrak g}

\newcommand\m{\mathfrak m}
\newcommand\n{\mathfrak n}
\renewcommand\u{\mathfrak u}

\newcommand\bC{\mathbb C}
\newcommand\bF{\mathbb F}

\newcommand\bP{\mathbb P}
\newcommand\bZ{\mathbb Z}

\newcommand\sub{\subseteq}

\renewcommand\bar{\overline}

\DeclareMathOperator{\Top}{top} %
\DeclareMathOperator{\GL}{GL} %
\DeclareMathOperator{\Mat}{Mat} %
\DeclareMathOperator{\Lie}{Lie} %
\DeclareMathOperator{\U}{U} %

\usepackage{fullpage}
\usepackage{amssymb}

\numberwithin{equation}{section}

\newtheorem{thm}{Theorem}[section]

\newtheorem{conj}[thm]{Conjecture}

\theoremstyle{definition}

\theoremstyle{remark}

\title[Calculating conjugacy classes in Sylow $p$-subgroups]
{Calculating conjugacy classes in Sylow $p$-subgroups \\ of finite
Chevalley groups}

\dedicatory{Dedicated to Gus Lehrer on the occasion of his 60th
birthday}

\author[S.~M.~Goodwin]{Simon M.~Goodwin}
\address{School of Mathematics, University of Birmingham, Birmingham, B15
2TT, United Kingdom}  \email{goodwin@maths.bham.ac.uk}
\urladdr{http://web.mat.bham.ac.uk/S.M.Goodwin}

\author[G.\ R\"ohrle]{Gerhard R\"ohrle}
\address{Fakult\"at f\"ur Mathematik, Ruhr-Universit\"at Bochum,
D-44780 Bochum, Germany} \email{gerhard.roehrle@rub.de}
\urladdr{http://www.ruhr-uni-bochum.de/ffm/Lehrstuehle/Lehrstuhl-VI/rubroehrle.html}

 


\begin{document}

\begin{abstract}
In \cite[\S 8]{Gconj}, the first author outlined an algorithm for
calculating a parametrization of the conjugacy classes in a Sylow
$p$-subgroup $U(q)$ of a finite Chevalley group $G(q)$, valid when
$q$ is a power of a good prime for $G(q)$. In this paper we develop
this algorithm and discuss an implementation in the computer algebra
language {\sf GAP}. Using the resulting computer program we are able
to calculate the parametrization of the conjugacy classes in $U(q)$,
when $G(q)$ is of rank at most $6$. In these cases, we observe that
the number of conjugacy classes of $U(q)$ is given by a polynomial
in $q$ with integer coefficients.
\end{abstract}

\maketitle

\section{Introduction}

Let $\U_n(q)$ be the subgroup of $\GL_n(q)$ consisting of upper
unitriangular matrices.  A longstanding conjecture, attributed to
G.~Higman \cite{H} states that the number of conjugacy
classes of $\U_n(q)$ is given by a polynomial in $q$ with integer
coefficients. This has been verified by computer calculation for $n
\le 13$ in the work of A.~Vera-L\'opez and J.~M.~Arregi, see
\cite{VA}. This conjecture has generated a great deal of interest,
see for example  \cite{robinson} and \cite{thompson}.

The equivalent problem of counting the number of (complex)
irreducible characters of $\U_n(q)$ has also attracted a lot of
attention, see for example \cite{lehrer}, \cite{isaacs1} and
\cite{isaacs2}. Thanks to work of M.~Isaacs \cite{isaacs1}, the
degrees of the irreducible characters of $\U_n(q)$ are known to all
be powers $q^d$ of $q$ and all exponents $d$ occur for $0 \le d
\le\mu(n)$, where $\mu(n)$ is an explicit upper bound due to work by
G.I.~Lehrer \cite{lehrer}.  It was conjectured by Lehrer
\cite{lehrer} that the number of irreducible characters of $\U_n(q)$
of degree $q^d$ is a polynomial in $q$ with integer coefficients
only depending on $n$ and $d$; this conjecture clearly implies
Higman's conjecture.

\smallskip

It is natural to consider the analogue of Higman's conjecture for
other finite groups of Lie type.  Below we introduce some notation
in order to discuss this analogue for the case of finite Chevalley
groups.

Let $G$ be a split simple algebraic group defined over the finite
field $\bF_p$ of $p$ elements, and assume that $p$ is good for $G$.
For a power $q$ of $p$ we write $G(q)$ for the finite group of
$\bF_q$-rational points of $G$; this is a finite Chevalley group.
Let $U$ be a maximal unipotent subgroup of $G$ defined over $\bF_p$,
so that $U(q)$ is a Sylow $p$-subgroup of $G(q)$.

In this paper, we describe an algorithm that calculates a
parametrization of the conjugacy classes of $U(q)$. We have
implemented this algorithm in the computer algebra language {\sf
GAP} \cite{GAP}. The algorithm is based on the outline given by the
first author in \cite{Gconj}. The output of the computer program
allows one to calculate the number $k(U(q))$ of conjugacy classes of
$U(q)$. Using our computer program we have proved the following
theorem.

\begin{thm}
\label{thm:main} Let $G$ be a split simple algebraic group defined
over $\bF_p$, where $p$ is good for $G$. Let $U$ be a maximal
unipotent subgroup of $G$ defined over $\bF_p$. Let $q$ be a power
of $p$.  If the rank of $G$ is at most $6$, then the number of
conjugacy classes of $U(q)$ is given by a polynomial in $q$ with
integer coefficients (and the polynomial itself is independent of
$p$).
\end{thm}

We have explicitly calculated the polynomial $k(U(q))$ for $G(q)$ of
rank at most $5$; these polynomials are presented in Table
\ref{Tab:kuq} in Section \ref{S:results}. For $G(q)$ of rank $6$,
the output of the computer program describing the parametrization of
the conjugacy classes of $U(q)$ is long and rather complicated.  It
is possible to check from this output in a short amount of time that
$k(U(q))$ is given by a polynomial in $q$; however, a great deal of
calculation would be required to explicitly compute this polynomial.
We explain this in more detail in Section \ref{S:results}.

Our computer program makes calculations using a $\bZ$-form of the
Lie algebra $\u$ of $U$.  In these calculations ``implicit
divisions'' are made, which lead to the output not being valid for a
finite number of primes $p$ that are recorded within the program. In
all case that we have run the program the only primes that are
output are bad primes for $G$.  A simple modification could insist
on not dividing by certain prime numbers, and in this way it is
possible to deal with the situation if any good primes were output.
We note that the results on which our algorithm is based require the
assumption that $p$ is a good prime for $G$, so the output of our
algorithm is not valid for bad primes even if there have been no
implicit divisions.

\smallskip

We have adapted the computer program so that it is also possible to
use it to calculate the number of $U(q)$-conjugacy classes in
certain subquotients of $U(q)$.  More precisely, let $B = N_G(U)$,
be the Borel subgroup of $G$ corresponding to $U$, then for normal
subgroups $N \supseteq M$ of $B$ that are contained in $U$, we can
calculate a parametrization of the $U(q)$-conjugacy classes in the
quotient $N(q)/M(q)$. We have made a number of such calculations in
case $G$ is of rank greater than $6$.  In all cases where we have
calculated the number of $U(q)$-conjugacy classes in $N(q)/M(q)$, we
observe that it is given by a polynomial in $q$.  In Table
\ref{Tab:kuql} in Section \ref{S:results}, we give the number
$k(U(q),U^{(l)}(q))$ of $U(q)$-conjugacy classes in the $l$th term
$U^{(l)}(q)$ of the descending central series of $U(q)$
for $G$ of exceptional type and certain $l$.

\smallskip

Generalizing a theorem of J.~Alperin \cite{A}, the authors showed in
\cite[Thm.\ 4.6]{GR} that if the centre of $G$ is connected, then
the number $k(U(q),G(q))$ of conjugacy classes of $U(q)$ in all of
$G(q)$ is a polynomial in $q$ with integer coefficients (in case
$G(q)$ has a simple component of type $E_8$, we require two
polynomials depending on the congruence of $q$ modulo $3$). The
theorem of Alperin \cite{A} can be viewed as support for Higman's
conjecture. Analogously, \cite[Thm.\ 4.6]{GR} suggests that, for $G$
not of type $E_8$, the number of conjugacy classes of $U(q)$ is
given by a polynomial in $q$. The results of our computer
calculations give supporting evidence for this behaviour, and we
thus propose the following analogue and extension of Higman's
conjecture for arbitrary finite Chevalley groups.

\begin{conj}
\label{conj:main}
Let $G$ be a split simple algebraic group defined
over $\bF_p$, where $p$ is good for $G$. Let $U$ be a maximal
unipotent subgroup of $G$ defined over $\bF_p$.
Let $q$ be a power of $p$.
If $G$ is not of type $E_8$, then $k(U(q))$ is
given by a polynomial in $q$ with integer coefficients (and the
polynomial itself is independent of $p$).
If $G$ is of type $E_8$, then $k(U(q))$ is given by
one of two polynomials depending on the congruence class of $q$ modulo $3$.
\end{conj}

The dependence of $k(U(q))$ on the congruence class of $q$ modulo
$3$ in the $E_8$ case in Conjecture \ref{conj:main} is suggested by
the $E_8$ case in \cite[Thm.\ 4.5(ii)]{GR}; though we do not wish to
rule out the possibility that there is just one polynomial. As
indicated above, due to the complexity of the computation, it is not
feasible to run our computer program in case $G$ is of type $E_8$.
In fact, as shown in Table \ref{Tab:kuql}, at present we have only
been able to calculate $k(U(q),U^{(l)}(q))$ explicitly for $l \ge
10$; we have $\dim U^{(10)} = 52$ and $\dim U = 120$ demonstrating
the difficulty of running our program for $E_8$.

\smallskip

As mentioned above our algorithm is not valid for bad primes.  It is
possible to calculate the $U(q)$-conjugacy classes for $G$ of type
$B_2$ and $p=2$ by hand.  In this case we have $k(U(q)) =
5(q-1)^2+4(q-1)+1$, which is a different polynomial to the one for
good primes given in Table \ref{Tab:kuq}; this is due to
degeneracies in the Chevalley commutator relations. In addition,
\cite[Thm.\ 4.6]{GR} is only valid for good primes, so we choose not
to make any conjecture for bad primes.

\smallskip From our calculations we can observe that each polynomial $k(U(q))$,
for $G$ of rank 5 or less, when written as a polynomial in $q-1$ has
non-negative integer coefficients, see Table \ref{Tab:kuq}. For $G$
of type $A_r$ and $r \le 12$ this was already observed in the
explicit results of Vera-L\'opez--Arregi  \cite{VA}. It would be
interesting to have a geometric interpretation of this positivity
behaviour. In Section \ref{S:results}, we give a reason why these
positivity phenomena hold for the cases that we have calculated.
This is done by analyzing the calculations made by the computer
program. We expect that if Conjecture \ref{conj:main} is true, then
the coefficients in $k(U(q))$ when written as a polynomial in $q-1$
are always non-negative.

In the cases where we have calculated $k(U(q))$, we have observed that
$k(U(q))$ always has constant term equal to $1$
when written as a polynomial in $q-1$.
In Section \ref{S:results}, we
explain why this is necessarily the case whenever $k(U(q))$ is a
polynomial in $q$.

Another observation is that the polynomial $k(U(q))$ is the same for
$G$ of type $B_r$ and $C_r$, for $r = 3,4,5$. It is likely that this
is always the case for any $r$.  We expect that this should be
explained by the duality of the underlying root systems.

Although $k(U(q),G(q))$ is always a polynomial in $q$ (assuming that
$G$ has connected centre and taking into account that for $G$ of
type $E_8$ there are two polynomials) and our calculations here tell
us that $k(U(q))$ is a polynomial for low rank, the number
$k(U(q),B(q))$ of $U(q)$-conjugacy classes in $B(q)$ is not always a
polynomial in $q$.  Indeed for $G$ of type $G_2$ it is shown in
\cite[Exmp.\ 4.8]{Gzeta} that the number $k(B(q),U(q))$ of
$B(q)$-conjugacy classes in $U(q)$ is given by two polynomials:
$q+15$ if $q$ is congruent to $-1$ mod $3$; and $q+17$ if $q$ is
congruent to $1$ mod $3$ (we are assuming $p$ is good for $G$).  A
general argument considering the commuting variety $\mathcal C(B,U)
= \{(b,u) \in B \times U \mid bu = ub\}$ shows that we always have
$k(U(q),B(q)) = (q-1)^2 k(B(q),U(q))$. Thus $k(U(q),B(q))$ is not a
polynomial in $q$ for $G$ of type $G_2$.

\medskip

Our algorithm calculates a family of varieties $X_c$ that
parameterize the conjugacy classes of $U$; moreover, these varieties
are defined over $\bF_p$.  The algorithm determines the polynomials
defining the $X_c$ as locally closed subsets of $(\bar\bF_p^\times)^{m_c}$ for
certain $m_c \in \bZ_{\ge 0}$.  The varieties $X_c$ are determined
with a backtrack algorithm using a depth-first search.  The
conjugacy classes of $U(q)$ are parameterized by the
$\bF_q$-rational points of the varieties $X_c$ and it is possible to
count these points.

The idea behind the algorithm is similar to that for the algorithm
used by B\"urgstein and Hesselink in \cite{BH} for calculating the
adjoint orbits of $B$ in $\u = \Lie U$; we remark that the algorithm
in \emph{loc.\ cit.\ }was not written to give a complete description
of the $B$-orbits in $\u$.  In addition, our algorithm generalizes
the one used in the work of Vera-L\'opez and Arregi for the type $A$
situation, see for example \cite{VA}.  Finally, we remark that the
algorithm of this paper uses ideas from the computer program
described in \cite{GDOOBS} in previous work of the first author.

\smallskip

We now give a brief outline of the structure of this paper.  In
Section \ref{S:recoll}, we introduce the notation that we require
and recall the relevant results from \cite{Gconj} and \cite{Gzeta}.
Then in Section \ref{S:alg} we describe the algorithm and its
implementation in {\sf GAP}. Finally, in Section \ref{S:results} we
discuss the results of our computations.  In particular, we present
explicit values for $k(U(q))$ for $G$ of rank at most $5$ (Table
\ref{Tab:kuq}) and the values of $k(U(q),U^{(l)}(q))$ for some cases
where $G$ is of exceptional type (Table \ref{Tab:kuql}).

\smallskip

As general references for algebraic groups defined over finite
fields we refer the reader to the books by Carter \cite{C} and
Digne--Michel \cite{DM}.

\section{Notation and recollection}
\label{S:recoll}

Let $p$ be a prime and let $G$ be a split simple algebraic group
defined over the finite field of $p$ elements $\bF_p$. We assume
throughout that $p$ is good for $G$ and we write $k$ for the
algebraic closure of $\bF_p$.

Fix a split maximal torus $T$ of $G$ and let $\Phi$ be the root
system of $G$ with respect to $T$. For a root $\alpha \in \Phi$ we
choose a generator $e_\alpha$ for the corresponding root subspace
$\g_\alpha$ of $\g = \Lie G$. Let $B \supseteq T$ be a Borel
subgroup of $G$ that is defined over $\bF_p$.  Let $U$ be the
unipotent radical of $B$ and let $\u = \Lie U$. Let $\Phi^+$ be the
system of positive roots of $\Phi$ determined by $B$.  The partial
order on $\Phi$ determined by $\Phi^+$ is denoted by $\preceq$.

For a power $q$ of $p$, we write $G(q)$ and $U(q)$ for the finite
groups of $\bF_q$-rational points of $G$ and $U$ respectively.  We
write $\u(q)$ for the Lie algebra of $\bF_q$-rational points of
$\u$.

\medskip

We now recall some results from \cite{Gconj} and \cite{Gzeta} on
which our algorithm for calculating the conjugacy classes of $U(q)$
is based.  Thanks to \cite[Thm.\ 1.1]{Gspring}, there are
generalizations of some results in \cite{Gconj}, as explained in
\cite[\S 6]{Gspring}; below we state the general versions without
further comment.  For the remainder of this section we fix a power
$q$ of $p$.

\smallskip

The theory of Springer isomorphisms can be used to show that the
conjugacy classes of $U(q)$ are in bijective correspondence with the
adjoint $U(q)$-orbits in $\u(q)$, see for example \cite[Prop.\
6.2]{Gconj}.  For the remainder of this paper, we will consider the
adjoint orbits of $U$ in $\u$ rather than the conjugacy classes of
$U$, as this is more convenient for our purposes.

\smallskip

Next we discuss the notion of minimal representatives of $U$-orbits
in $\u$, and how they are used to partition the set of $U$-orbits in
$\u$.  The reader is referred to \cite[\S 5 and \S 6]{Gconj} and
\cite[\S 3 and \S 4]{Gzeta} for full details.

We fix an enumeration of the set of positive roots $\Phi^+ =
\{\beta_1,\dots,\beta_N\}$, such that $i \le j$ whenever $\beta_i
\preceq \beta_j$, and define the sequence of $B$-submodules
$$
\u = \m_0 \supseteq \dots \supseteq \m_N = \{0\}
$$
of $\u$ by $\m_i = \sum_{j=i+1}^N \g_{\beta_j}$.  We consider the
action of $U$ on successive quotients $\u_i = \u/\m_i$ induced from
the adjoint action of $U$ on $\u$. We note that the parametrization
of the adjoint $U$-orbits described below depends on the choice of
the enumeration of $\Phi^+$.

Let $x \in \u$ and consider the set
$$
x + ke_{\beta_i} + \m_i = \{x + \lambda e_{\beta_i} + \m_i \mid
\lambda \in k \} \sub \u_i.
$$
By \cite[Lem.\ 5.1]{Gconj}, for $x \in \u$ either:
\begin{itemize}
\item[(I)] all elements of $x + ke_{\beta_i} + \m_i$ are $U$-conjugate; or
\item[(R)] no two elements of $x + ke_{\beta_i} + \m_i$ are $U$-conjugate.
\end{itemize}
We say that
\begin{itemize}
\item $i$ is an {\em inert point} of $x$ if (I) holds;
\item $i$ is a {\em ramification point} of $x$ if (R) holds.
\end{itemize}
An element $x + \m_i= \sum_{j=1}^i x_j e_{\beta_j}+ \m_i$ of $\u_i$
is said to be the {\em minimal representative} of its $U$-orbit in
$\u_i$ if $x_j=0$ whenever $j$ is an inert point of $x$.  It follows
from \cite[Prop.\ 5.4 and Lem.\ 5.5]{Gconj} that each $U$-orbit in
$\u_i$ contains a unique minimal representative; in particular, this
holds for the action of $U$ on $\u$.

Thanks to \cite[Prop.\ 4.2 and Lem.\ 5.7]{Gconj}, we have that $i$
is an inert point of $x \in \u$ if and only if $\dim \c_\u(x+\m_i) =
\dim \c_\u(x+\m_{i-1}) - 1$; if $i$ is a ramification point of $x$,
then we have $\dim \c_\u(x+\m_i) = \dim \c_\u(x+\m_{i-1})$.  Here
$\c_\u(x+\m_i)$ is the centralizer of $x + \m_i$ for the action of
$\u$ on $\u_i$ induced from the adjoint action of $\u$ on itself.

The above discussion implies that the adjoint orbits of $U$ in $\u$
are parameterized by their minimal representatives.  Further, the
set of minimal representatives can be partitioned into sets $X_c$
for $c \in \{\I,\R\}^N$: the set $X_c$ is defined to consist of the
minimal representatives $x \in \u$ of the $U$-orbits in $\u$ such
that for all $i = 1,\dots,N$ we have that $i$ is an inert point of
$x$ if and only if $c_i = \I$.  Thanks to \cite[Prop.\ 2.4]{Gzeta},
each of the sets $X_c$ is a locally closed subset of $\u$, and
therefore has the structure of an algebraic variety.

The above partition of the $U$-orbits in $\u$ can be refined to be
indexed by $N$-tuples $c \in \{\I,\R_\nn,\R_0\}^N$ as follows.  For
$c \in \{\I,\R_\nn,\R_0\}^N$, the set $X_c$ is defined to consist of
the minimal representatives $x = \sum x_i e_{\beta_i} \in \u$ of the
$U$-orbits in $\u$ such that for all $i = 1,\dots,N$ we have that
$i$ is an inert point of $x$ if and only if $c_i = \I$; and if $c_i
\ne \I$, then $x_i = 0$ if and only if $c_i = \R_0$. Thanks to
\cite[Lem.\ 4.2]{Gzeta}, each of the sets $X_c$ is a locally closed
subset of $\u$, and therefore has the structure of an algebraic
variety.  In fact, $X_c$ is a subvariety of $\{(x_j)_{c_j = \R_\nn}
\mid x_j \in k^\times\} \cong (k^\times)^{m_c}$, where $m_c = |\{j
\mid c_j = \R_\nn\}|$.

\medskip

We now explain how the above parametrization of the $U$-orbits in $\u$
descends to give a parametrization of the $U(q)$-orbits in $\u(q)$.  The
reader is referred to \cite[\S 6]{Gconj} for further details.

Thanks to \cite[Prop.\ 4.5]{Gconj}, we have that for $x \in \u$ the
centralizer $C_U(x)$ of $x$ in $U$ is connected.  This implies that
the $U(q)$-orbits in $\u(q)$ correspond bijectively to the
$U$-orbits in $\u$ that are defined over $\bF_q$.  Let $x \in \u$ be
the minimal representative of its $U$-orbit. Then, by \cite[Lem.\
6.3]{Gconj}, the orbit $U \cdot x$ is defined over $\bF_q$ if and
only if $x \in \u(q)$.  We require that the definition of $G$ over
$\bF_q$ is split for this last assertion, and this is the reason for
this assumption.  If the definition of $G$ over $\bF_q$ is not
split, then it is a non-trivial task to determine which minimal
representatives correspond to orbits defined over $\bF_q$.

It follows from the above discussion that the adjoint orbits of
$U(q)$ in $\u(q)$ are parameterized by the minimal representatives
of the $U$-orbits in $\u$ that lie in $\u(q)$.  In turn these
minimal representatives are partitioned into the sets $X_c(q)$ of
$\bF_q$-rational points of the varieties $X_c$, for $c \in
\{\I,\R_\nn,\R_0\}^N$.

\section{The algorithm}
\label{S:alg}

In this section we develop the algorithm outlined in \cite{Gconj}
for calculating the parametrization of the adjoint $U$-orbits in
$\u$. The idea is to calculate the polynomials defining the
varieties $X_c$ for $c \in \{\I,\R_\nn,\R_0\}^N$ as locally closed
subsets of $(k^\times)^{m_c}$.  We present the algorithm, and then
explain why the algorithm does indeed calculate a parameterization
of the adjoint $U$-orbits in $\u$.  Next we discuss two
modifications that are used in the {\sf GAP} implementation, before
briefly explaining the implementation. Finally, we explain how the
output of the computer program is used to calculate $k(U(q))$.

\smallskip

In order to explain the algorithm we have to introduce some more
notation; we continue to use the notation given in the previous
section.

We wish to consider all primes $p$ simultaneously, so we need a
$\bZ$-form of $\g$.  Let $\g_\bC$ be the complex simple Lie algebra
of the same type as $\g$.  Fix a Chevalley basis of $\g_\bC$ and let
$\g_\bZ$ be the corresponding $\bZ$-form of $\g_\bC$.  We let
\[
m := \max \{ m_c \mid c \in \{\I,\R_\nn,\R_0\}^N, X_c \ne
\varnothing\},
\]
where $m_c = |\{i \mid c_i = \R_\nn\}|$, as defined earlier.   We
define
$$
\tilde \u = \u_\bZ \otimes_\bZ \bZ[t_1,\dots,t_m],
$$
where $\bZ[t_1,\dots,t_m]$ is the polynomial ring in $m$
indeterminates $t_1,\dots,t_m$.  We denote by
$e_{\beta_1},\dots,e_{\beta_N}$ the elements of the Chevalley basis
of $\g_\bZ$ corresponding to $\Phi^+$, which is enumerated as in the
previous section. These elements form a $\bZ$-basis of $\u_\bZ$, and
by a minor abuse we also consider them as elements of both $\u$ and
$\tilde \u$.

Let $c \in \{\I,\R_\nn,\R_0\}^i$ for some $i \le N$. For $j =
1,\dots,m_c$, we define $\beta_{c,j} \in \Phi^+$ by setting
$\beta_{c,j} = \beta_l$, where $l$ is the $j$th smallest element in
$\{ m \mid c_m = \R_\nn\}$.  We associate to $c$ the element
$$
x_c(t) = \sum_{j=1}^{m_c} t_j e_{\beta_{c,j}} \in \tilde \u.
$$
Given $\tau = (\tau_1,\dots,\tau_{m_c}) \in k^{m_c}$, we write
$x_c(\tau)$ for the element of $\u$ obtained by substituting $t_j =
\tau_j$ in $x_c(t)$, i.e.
$$
x_c(\tau) = \sum_{j=1}^{m_c} \tau_j e_{\beta_{c,j}}.
$$

The variety $X_c$ is a locally closed subset of $\{x_c(\tau) \mid
\tau \in (k^\times)^{m_c}\} \cong (k^\times)^{m_c}$.  Therefore,
there are subsets $A_c^1,\dots,A_c^{l_c}$ and
$B_c^1,\dots,B_c^{l_c}$ of $k[t_1,\dots,t_{m_c}] \sub
k[t_1,\dots,t_m]$ such that $X_c$ is the disjoint union of the sets
$$
X_c^i = \{x_c(\tau) \mid f(\tau) = 0 \text{ for all } f \in A_c^i
\text{ and } g(\tau) \ne 0 \text{ for all } g \in B_c^i\},
$$
for $i=1,\dots,l_c$. In fact the polynomials in the sets $A_c^i$ and
$B_c^i$ can be taken to have integer coefficients; this is due to
the integrality of the Chevalley commutator relations.  The purpose
of our algorithm is to determine certain choices for the sets
$A_c^i$ and $B_c^i$.

We note here that it is most often the case that we can take $l_c =
1$ and $A_c^1 = B_c^1 = \varnothing$.  The values of $c$ for which
this is not the case in some sense explain why the determination of
the conjugacy classes of $U(q)$ is complicated in general.  We also
remark that it is often the case that $X_c = \varnothing$, which
corresponds to the case $l = 0$.

\smallskip

We now introduce some notation needed in order to say how the sets
$A_c^i$ and $B_c^i$ are determined in the algorithm. Let
$y_1,\dots,y_N$ be variables. We may write:
$$
\left[\sum_{j=1}^N y_j e_{\beta_j},x_c(t)\right] = \sum_{j=1}^N
\sum_{k=1}^N P_{jk}^c(t) y_k e_{\beta_j},
$$
where each $P_{jk}^c(t) \in \bZ[t_1,\dots,t_m]$ is linear: this is
easily achieved using the Chevalley commutator relations for
$\g_\bC$. It is then the case that $\dim \c_\u(x_c(\tau)+\m_i)$ is
the dimension of the solution space of the system of linear
equations:
$$
\sum_{k=1}^N P_{jk}^c(\tau) y_k = 0,
$$
for $j=1,\dots,i$.

\smallskip

We are now in a position to describe our algorithm.  It calculates
sets of polynomials that determine certain choices of the varieties
$X_c^i$.  These polynomials are calculated using the $\bZ$-form
$\u_\bZ$ of $\u$ and they have rational coefficients.
Observe that in step (2)(a) of the algorithm
there may be ``implicit divisions'' by certain primes for which the
output will not be valid.  This is discussed in more detail later in
this section.

The algorithm uses a backtrack algorithm with a depth-first search
to calculate certain choices of sets $A_c^i$ and $B_c^i$ that
determine the varieties $X_c^i$ for each $c \in
\{\I,\R_\nn,\R_0\}^N$, and $i = 1,\dots,l_c$. In the algorithm we
require a total order on $\bZ[t_1,\dots,t_m]$; we use the order
defined in precedence by the number of terms, total degree and
leading coefficient (with respect to the degree then lexicographic
order on monomials).

The variables used in the algorithm are:
\begin{itemize}
\item the ``current string'' $c$ is an element of $\{\I,\R_\nn,\R_0\}^i$
for some $i$ and determines $x_c(t) \in \tilde \u$ as above;
\item the set of ``satisfied'' polynomials $A$ is a
subset of $\bZ[t_1,\dots,t_m]$;
\item the set of ``unsatisfied''
polynomials $B$ is a subset of $\bZ[t_1,\dots,t_m]$;
\item the matrix $Q(t)$ is an element of $\Mat_{i,N}(\bZ[t_1,\dots,t_m])$,
which is obtained from the matrix $(P_{jk}^c(t))$ by ``row reducing'' the first
$i$ rows;
\item the ``pivot string'' $\pi$ is an element of $\{0,1,\dots,N\}^i$, which
``records the columns used in the row reductions'';
\item the stack $S = \{(c, A, B, \pi, Q(t))\}$ is an (ordered) subset of
$$
\bigcup_{i=1}^N \{\I,\R_\nn,\R_0\}^i \times
\bP(\bZ[t_1,\dots,t_m])^2 \times \bigcup_{i=1}^N \{0,1,\dots,N\}^i
\times \Mat_N(\bZ[t_1,\dots,t_m]),
$$
which contains variables to be considered later in the algorithm;
and
\item the output set $O$ is a
subset of $\{\I,\R_\nn,\R_0\}^N \times \bP(\bZ[t_1,\dots,t_m])^2$.
\end{itemize}
Here, $\bP$ denotes the power set.

The stack is required to be ordered as the algorithm takes
elements from the ``top'' of the stack. The element at the top is
the one that has been most recently added and is denoted by
$\Top(S)$. The initial configuration in the algorithm is
as follows:
\begin{itemize}
\item $c := (\R_\nn)$;
\item $A := \varnothing$;
\item $B := \varnothing$;
\item $\pi : = (0)$;
\item $Q(t)$ is the $1 \times N$ matrix with all entries equal to $0$;
\item $S := \{(\R_0,\varnothing,\varnothing,(0),Q(t))\}$; and
\item $O := \varnothing$.
\end{itemize}

Now we explain the next step in the algorithm; we have numbered the
steps in the algorithm, so that we can refer back to it in the
explanation given afterwards.

\begin{itemize}
\item[(1)] If the length of $c$ is $N$, then we are finished with this
string.  We set:
\begin{itemize}
\item[(a)] $O := O \cup \{(c,A,B)\}$
\end{itemize}
If $S = \varnothing$, then we finish.
\\
Else we make the following changes to the variables:
\begin{itemize}
\item[(b)] $(c,A,B,\pi,Q) := \Top(S)$; and
\item[(c)] $S := S \setminus \{\Top(S)\}$.
\end{itemize}
\item[(2)]  If the length of $c$ is $i-1 < N$, then we proceed by making
the $i$th row reduction for the matrix $(P_{jk}^c(t))$ as defined
above. Note that $Q(t)$ is the matrix resulting from the first $i-1$
row reductions.  We first append the $i$th row of $(P_{jk}^c(t))$ to
$Q(t)$ and then make the row reduction as follows:
\begin{itemize}
\item[(a)] for $j=1,\dots,i-1$, if $\pi_j \ne 0$
we set $Q_i(t):=Q'_{j,\pi_j}(t)Q_i(t) - Q'_{i,\pi_j}(t)Q_j(t)$, where
$Q'_{j,\pi_j}(t)$ is $Q_{j,\pi_j}(t)$ divided by the highest common factor
of $Q_{j,\pi_j}(t)$ and $Q_{i,\pi_j}(t)$, and $Q'_{i,\pi_j}(t)$ is defined
analogously.
\end{itemize}
Let $R_i$ be the set of non-zero polynomials in the $i$th row of
$Q(t)$ that are not divisible by any element of $A$. We next
consider three cases:
\begin{itemize}
\item[(b)] $R_i = \varnothing$. We update the variables as
follows:
\begin{itemize}
\item[(i)] $\pi:=(\pi,0)$
\item[(ii)] $c:=(c,\R_\nn)$; and
\item[(iii)] $S:=S \cup \{((c,\R_0),A,B,\pi,Q)\}$.
\end{itemize}
\item[(c)] $R_i \ne \varnothing$ and there is some non-zero
element of $R_i$ that is a monomial or divides some element of $B$.
We let $Q_{i,l}(t)$ be the least such polynomial with respect to our
chosen order on the set of all polynomials.  We update the variables
as follows:
\begin{itemize}
\item[(iv)] $c:=(c,\I)$; and
\item[(v)] $\pi:=(\pi,l)$.
\end{itemize}
\item[(d)]  Otherwise, we pick a least element $Q_{i,l}(t)$ of $R_i$. We update the
variables as follows:
\begin{itemize}
\item[(vi)] $S:=S \cup \{(c,A \cup \{Q_{i,l}(t)\},B,\pi,Q)\}$;
\item[(vii)] $c:=(c,\I)$;
\item[(viii)] $B:=B \cup \{Q_{i,l}(t)\}$; and
\item[(ix)] $\pi := (\pi,l)$.
\end{itemize}
\end{itemize}
\end{itemize}

The output of the algorithm is a collection of triples $(c,A,B)$.
Each of these triples determines a subvariety
\begin{equation} \label{e:X_c,A,B}
X_{c,A,B} = \left\{\sum_{j : c_j = \R_\nn} \tau_j e_{\beta_{c,j}}
\in X_c \mid f(\tau) = 0 \text{ for all } f \in A \text{ and }
g(\tau) \ne 0 \text{ for all } g \in B\right\}
\end{equation}
of $X_c$.  For fixed $c$, we can write $A_c^1,\dots,A_c^{l_c}$ and
$B_c^1,\dots,B_c^{l_c}$ for the sets $A$ and $B$ occurring in a
triple $(c,A,B)$.  These are the determined choice of the sets
$A_c^i$ and $B_c^i$ such that $X_c$ is the disjoint union of the
varieties $X_c^i$ as defined earlier; each $X_c^i$ is equal to the
corresponding $X_{c,A,B}$.

We now explain why the output of our algorithm can be used to
determine all the minimal representatives of the $U(q)$-orbits in
$\u(q)$, for almost all primes $p$.  As explained above there may be
``implicit divisions'' in step (2)(a) of the algorithm leading to a
finite number of primes for which we cannot determine the minimal
representatives of the $U(q)$-orbits in $\u(q)$.  We choose not to
give a formal proof of the correctness of the algorithm as this
would be very technical and just give an outline.  We argue by
induction on $i$ to show that the algorithm determines varieties
$X_{c,A,B,i}$ (as defined below) from which one can calculate all
the minimal representatives of the $U(q)$-orbits in $\u_i(q)$ for
each $i$ (and valid $p$).  We do not discuss the part of the
algorithm dealing with the row reductions in (2)(a), as this is
elementary.

In the discussion below we often speak of all \emph{relevant} $\tau
\in k^m$. When considering the triple $(c,A,B)$, this means all
$\tau \in k^m$ such that $f(\tau) = 0$ for all $f \in A$ and
$g(\tau) \ne 0$ for all $g \in B$, i.e.\ all $\tau \in k^m$ for
which $x_c(\tau) \in X_{c,A,B}$.

We need to explain our inductive hypothesis as this is not
compatible with the depth-first search used in the algorithm.  To do
this we must consider all triples $(c,A,B)$ with $c$ of length $i$
that occur during the running of the algorithm.  For such $(c,A,B)$
varieties, $X_{c,A,B}$ can be defined as in \eqref{e:X_c,A,B}.  Then
the inductive hypothesis says that the varieties $X_{c,A,B,i} = \{x
+ \m_i \mid x \in X_{c,A,B}\}$ give all minimal representatives of
the $U$-orbits in $\u_i$.

So assume the inductive hypothesis for $i-1$.  Then we have to show
that for each $(c,A,B)$ with $c$ of length $i$ and all minimal
representatives of $U$-orbits in $\u_i$ of the form $x_c(\tau) +
\lambda e_{\beta_i}+ \m_i$ with $x_c(\tau) \in X_{c,A,B}$ and
$\lambda \in k$, we have that $x_c(\tau) + \lambda e_{\beta_i}$ lies
in some variety of the form $X_{c',A',B'}$ where $c' = (c,Z)$ and $Z
\in \{\I,\R_0,\R_\nn\}$. To do this we have to consider the steps in
(2) in the algorithm.

After the row reduction made in (2)(a), we have the row reduced
matrix $Q(t)$ and the set $R_i$. As explained earlier we have that
$i$ is an inert point of $x_c(\tau) \in X_{c,A,B}$ if and only if
$\dim \c_u(x_c(\tau)+\m_i) = \dim \c_u(x_c(\tau)+\m_{i-1}) - 1$.
Also the dimension of $\c_u(x_c(\tau)+\m_i)$ is the rank of the row
reduced matrix $Q(\tau)$. Therefore, $i$ is an inert point of
$x_c(\tau)$ if and only if $f(\tau) \ne 0$ for some $f \in R_i$.
Note that we only have to consider the polynomials in $R_i$, because
for any non-zero entry in the $i$th row of $Q(t)$ that is divisible
by some polynomial $f(y) \in B$, we automatically have $f(\tau)=0$
for all relevant $\tau$. At this stage one ideally wants to
determine for which values of $\tau$ there is some $f \in R_i$ with
$f(\tau) \ne 0$. However, this is a difficult task if $R_i$ contains
several polynomials. Also we require a fixed value for the next
entry of $\pi$, so that we can perform the row reductions later in
the algorithm. So we proceed by considering the three cases
(2)(b)--(d) in the algorithm:
\begin{enumerate}
\item[(b)]  In this case it is clear that $i$ is a ramification point
of $x_c(\tau)$ for all relevant $\tau$. Therefore, $x_c(\tau) +
\lambda e_{\beta_i}+\m_i$ is the minimal representative of its
$U$-orbit in $\u_i$ for all $\lambda \in k$. We see that the strings
$c'$ corresponding to all such minimal representatives are passed on
in the program: the non-zero values of $\lambda$ correspond to the
updated strings $c' = (c,\R_\nn)$ in (ii), and the case $\lambda =
0$ corresponds to the string $c' = (c,\R_0)$ added to the stack.
\item[(c)]  For this case, it is clear that $f(\tau) \ne 0$ for all
relevant $\tau$, where $f(t) = Q_{i,l}(t)$ is the least element of
$R_i$ that is either a monomial or divides some element of $B$.
Therefore, $i$ is an inert point of $x_c(\tau)$ for all relevant
$\tau$.  Thus the only minimal representative of its $U$-orbit in
$\u_i$ of the form $x_c(\tau) + \lambda e_{\beta_i} + \m_i$ is when
$\lambda = 0$. This minimal representative corresponds to the
updated string $c' = (c,\I)$ in (iv).
\item[(d)]  This case is more complicated.  We consider the least element
$f(t) = Q_{i,l}(t)$ of $R_i$ and the following cases for $\tau$:
\begin{enumerate}
\item[(I)] If $f(\tau) \ne 0$, then $i$ is an inert point for $x_c(\tau)$.
For such $\tau$, the only minimal representative of its $U$-orbit in
$\u_i$ of the form $x_c(\tau) + \lambda e_{\beta_i} + \m_i$ is when
$\lambda = 0$.  This minimal representative corresponds to the
updated string $c' = (c,\I)$ in (vii), along with $f$ being added to
$B$ in (viii).
\item[(II)]  If $f(\tau) = 0$, then we cannot say whether $i$ is an inert
or ramification point of $x_c(\tau)$.  The element
$(c,A\cup\{f\},B,\pi,Q(t))$ added to the stack in (vi) will be
considered later in the algorithm.  For this case $R_i$ will have
fewer elements, and will be considered again either in case (a) or
(c). This process will finish, and it is straightforward to see that
when this happens all the required triples $(c',A',B')$ with $c' =
(c,Z)$, will have occurred.
\end{enumerate}
\end{enumerate}
Putting together the above case analysis, verifies the inductive
step.  Therefore, for all minimal representatives $x_c(\tau)$ of the
$U$-orbits in $\u$, a corresponding triple $(c,A,B)$ occurs at some
point of the program.  This triple is added to the output in (1)(a),
which means that the output of the algorithm determines all minimal
representatives, as claimed.

\smallskip

Next we discuss two modifications to the algorithm that we make for
its implementation in {\sf GAP}.  These changes are made in order to
speed up the computations; we chose not to include them in the above
description of the algorithm for simplicity.

Our first modification allows us to reduce the complexity by using
the action of the maximal torus $T$ to ``normalize'' certain
coefficients to be equal to 1.  Let $c \in \{\I,\R_\nn,\R_0\}^N$ and
let $x_c(t) = \sum_{i=1}^{m_c} t_i e_{\beta_{c,i}}$ be defined as
above. Suppose that $\{\beta_{c,i} \mid i \in J \}$ is linearly
independent for some $J \sub \{1,\dots,m_c\}$. Then for every  $\tau
\in (k^\times)^{m_c}$ there is some $\sigma \in (k^\times)^{m_c}$
with $\sigma_i = 1$ for all $i \in J$ and $x_c(\tau)$ is conjugate
to $x_c(\sigma)$ via $T$, i.e.\ the action of $T$ can be used to
``normalize $\tau_i = 1$'' for all $i \in J$. In the computer
program we replace $x_c(t)$ by $\sum_{i \in J} e_{\beta_{c,i}} +
\sum_{i \notin J} t_i e_{\beta_{c,i}}$ thereby reducing the number
of indeterminates required.  This speeds up the row reduction of the
matrix $Q(t)$ significantly.   Some care is needed when using this
modification, as the centralizer $C_T(x_c(\tau))$ can be
disconnected.  If this is the case then it can become difficult to
determine the elements of $X_c(q)$ from the ``normalized'' elements
of $X_c$.  This problem can be resolved by not allowing certain
normalizations; we omit the technical details here.

The second adaptation deals with ``easy'' elements of the set $A$.
If there is a linear polynomial in $A$, then we may simplify future
checks by ``making a substitution''.  If $t_i - a(t) \in A$, where
$i \in \{1,\dots,m_c\}$ and $a(t)$ is a linear polynomial not
involving $t_i$, then we can make the substitution $t_i = a(t)$ in
$x_c(t)$,  in the polynomials in $A$ and $B$, and in the matrix
$Q(t)$; we then remove $t_i - a(t)$ from $A$.  This modification
reduces the number of indeterminates, and also the number of
polynomials in $A$.  This helps to speed up the program.

\smallskip

We next explain a check that has to be included in the program to
see which primes the output is valid for. When making the row
reductions in (2) of our algorithm, there are implicit divisions by
certain integers. Essentially, we need to be able to divide by the
polynomials in the set $B$ in (2)(a). Therefore, we keep track of
the primes dividing their leading coefficients, and these primes are
output by the program.  The output of our program cannot be used to
determine the minimal representatives of the $U(q)$-orbits in
$\u(q)$ for such $p$. The only primes output by the program for the
cases that we have calculated are bad primes for $G$.  If a good
prime $p$ were output by the program, then it would be
straightforward to adapt the program to insist that there are no
implicit divisions by $p$, and running this modification would give
an output valid for $p$.

\smallskip

The algorithm is implemented in {\sf GAP} with the two modifications
and the check for primes.  This is achieved using the functions for
Lie algebras and polynomial rings in {\sf GAP}. This allows us to
define $\tilde \u$ within {\sf GAP} and therefore allows us to
obtain the matrices $(P_{jk}^c(t))$ that we row reduce using the
method given in (2)(a) of the algorithm. The implementation is based
on the algorithm given and the two modifications discussed above. We
choose not to include any of the technical details.

\smallskip

In the next section we present the values of $k(U(q))$ that we have
calculated from the output of our program.  Each of the varieties
$X_{c,A,B}$ is defined by polynomials with integer coefficients, so
is defined over $\bF_p$. We have
$$
k(U(q)) = \sum_{(c,A,B)} |X_{c,A,B}(q)|.
$$
We can, therefore, calculate $k(U(q))$ by calculating
$|X_{c,A,B}(q)|$ for all triples $(c,A,B)$.  If the polynomials in
$A$ and $B$ are not too complicated, then this can be achieved quite
easily.  We discuss this below.

It is most commonly the case that $c$ occurs in just one triple
$(c,A,B)$ for which both $A$ and $B$ are empty. In which case it is
easily seen that $X_c = X_{c,A,B}$ and $|X_c(q)| = (q-1)^{m_c}$. The
next simplest case is when $A \cup B$ has one element that is
linear. For example, consider the polynomial $t_1-1$: if $A =
\{t_1-1\}$ and $B = \varnothing$, then $|X_{c,A,B}(q)| =
(q-1)^{m_c-1}$; and if $A = \varnothing$ and $B = \{t_1-1\}$, then
$|X_{c,A,B}(q)| = (q-1)^{m_c-1}(q-2)$.  More complicated sets $A$
and $B$ that we need to consider require a little thought to
calculate $|X_{c,A,B}(q)|$.

As the rank of $G$ increases the polynomials become more
complicated. For the $F_4$, $B_5$ and $C_5$ cases we get a number of
quadratic polynomials. For the rank $6$ cases, the polynomials
become more complicated still and the number of triples $(c,A,B)$
with $A$ or $B$ non-empty gets large.  From the output of the
program we can view all the polynomials occurring in the sets $A$
and $B$.  One can check that for each of these polynomials it is
possible to ``solve for one indeterminate in terms of the others''.
With a little further consideration one can see that this means each
of the sets $X_{c,A,B}(q)$ has size a polynomial in $q$.  However,
the number of such $X_{c,A,B}(q)$ is so large that it would be
rather time consuming to calculate $k(U(q))$ explicitly.

\section{Results}
\label{S:results}

In this final section we present some explicit results of our computations and
go on to discuss some interesting features of the output.

In Table \ref{Tab:kuq} below we present the polynomials $k(U(q))$
for $G(q)$ of rank at most $5$; in this table we let $v = q-1$ to
save space. We include the values for $G$ of type $A_r$ for
completeness, though these polynomials have been known for some
time, thanks to the work of Vera-L\'opez and Arregi referred to in the
introduction. Also, as discussed below, the value of $k(U(q))$ is
the same for $G$ of type $B_r$ and $C_r$, so we only include this
polynomial once.

\begin{table}[h!tb]
\renewcommand{\arraystretch}{1.5}
\begin{tabular}{|l|l|}
\hline
$G$ & $k(U(q))$  \\
\hline\hline
$A_1$ & $v+1$ \\
\hline
$A_2$ & $v^2+3v+1$ \\
$B_2$ & $2v^2+4v+1$ \\
$G_2$ & $v^3+5v^2+6v+1$ \\
\hline
$A_3$ & $2v^3+7v^2+6v+1$ \\
$B_3$/$C_3$ & $v^4+8v^3+16v^2+9v+1$ \\
\hline
$A_4$ & $5v^4+20v^3+25v^2+10v+1$\\
$B_4$/$C_4$ & $v^6+11v^5+48v^4+88v^3+64v^2+16v+1$ \\
$D_4$ & $2v^5+15v^4+36v^3+34v^2+12v+1$\\
$F_4$ & $
v^8+9v^7+40v^6+124v^5+256v^4+288v^3+140v^2+24v+1$\\
\hline
$A_5$ & $v^6+18v^5+70v^4+105v^3+65v^2+15v+1$\\
$B_5$/$C_5$ &
$2v^8+24v^7+132v^6+395v^5+630v^4+500v^3+180v^2+25v+1$ \\
$D_5$ & $2v^7+22v^6+106v^5+235v^4+240v^3+110v^2+20v+1$\\
\hline
\end{tabular}
\medskip
\caption{$k(U(q))$, as polynomials in $v = q-1$} \label{Tab:kuq}
\end{table}

We make some comments about the polynomials in Table \ref{Tab:kuq}.
We start by making the observation that $k(U(q))$ considered as a
polynomial in $v = q-1$ has non-negative coefficients. For the case
$G$ is of type $A_r$ ($r \le 12$), this was observed by Vera-L\'opez
and Arregi in \cite{VA}.  It would be interesting to have a
geometric explanation of these positivity phenomena.

We give a heuristic idea why this occurs for the cases that we have
calculated by considering the partition of the conjugacy classes
used in our algorithm. As discussed at the end of the previous
section, the number of the $\bF_q$-rational points of the varieties
$X_c^i$ is most commonly $|X_c^i(q)| = v^{m_c}$. Although there are
some values of $c$ and $i$ for which $|X_c^i(q)|$ is a polynomial in
$v$ with negative coefficients, these negative coefficients are few
enough so that they are cancelled by the families of size $v^{m_c}$.

\smallskip

We observe that the constant coefficient in $k(U(q))$ as a
polynomial in $v$ is always $1$.  This is explained by the action of
the split maximal torus $T$ of $G$ on each $X_c$ for all $c \in
\{\I,\R_\nn,\R_0\}^N$. This action is non-trivial unless $c_i =
\R_0$ for all $i$, so that $X_c = \{0\}$. It is easy to see that if
$X_c \ne \{0\}$, then the orbits of $T(q)$ on $X_c(q)$ are all of
size $v^a/b$ for some $a, b \in \bZ_{\ge 1}$, so $|X_c(q)|$ is
divisible by $v^a/b$.  This implies that the constant coefficient in
$k(U(q))$ as a polynomial in $v$ must be 1 (corresponding to the
zero orbit).

\smallskip

We now comment on the fact that the value of $k(U(q))$ is the same
for $G$ of type $B_r$ and $C_r$, for $r = 3,4,5$.  One can see that
the groups $U(q)$ are not isomorphic for $G$ of type $B_r$ and
$C_r$: thanks to a result of A.~Mal'cev \cite{malcev}, the maximal
size of an abelian subgroup of $U(q)$ is different for $G$ of type
$B_r$ and $C_r$. Using the variation of our program discussed below,
one can also show that the number of $U(q)$-conjugacy classes in the
derived subgroup $U^{(1)}(q)$ of $U(q)$ are different for $G$ of
types $B_r$ and $C_r$, for $r = 3,4,5$.  It would be interesting to
have a reason for the coincidences in the numbers $k(U(q))$; we
expect it should be explained by the duality of the root systems of
type $B_r$ and $C_r$, see for example \cite[Ch.\ 4]{C} for similar
phenomena.

\smallskip

As mentioned in the introduction, we have adapted our program to
consider the action of $U$ on certain subquotients $M/N$.  The
adaption is valid when $M \supseteq N$ are normal subgroups of $B$
contained in $U$.  The algorithm runs in essentially the same way:
one has to replace the filtration of $\u$ by an analogous filtration
of $\m/\n$, then change the initial configuration and the point at
which variables are added to the output set $O$ accordingly.

In Table \ref{Tab:kuql} below we give some values of
$k(U(q),U^{(l)}(q))$ for $G$ of exceptional type. We recall that the
descending central series of $U$ is defined by $U^{(0)} = U$ and
$U^{(l)} = [U^{(l-1)},U]$ for $l \ge 1$.  The cases that we have
included are those for which we are able to compute
$k(U(q),U^{(l)}(q))$ in a reasonable amount of time and for which
there is an infinite number of $B$-orbits in $\u^{(l)} = \Lie U^{(l)}$;
we refer the reader to \cite{GRfin} for a classification of all cases
when there is only a finite number of $B$-orbits in $\u^{(l)}$
for $G$ of exceptional type.

\begin{table}[h!tb]
\renewcommand{\arraystretch}{1.5}
\begin{tabular}{|l|r|l|}
\hline
$G$ & $l$ & $k(U(q),U^{(l)}(q))$  \\
\hline\hline
$F_4$ & $1$ & $v^7+7v^6+24v^5+63v^4+119v^3+88v^2+20v+1$ \\
& 2 & $2v^5+14v^4+50v^3+58v^2+17v+1$ \\
& 3 & $2v^4+18v^3+35v^2+14v+1$ \\
\hline
$E_6$ & 1 &
$v^{10}+10v^9+47v^8+153v^7+435v^6+993v^5+1315v^4+868v^3+255v^2+30v+1$\\
& 2 & $2v^7+28v^6+160v^5+386v^4+404v^3+165v^2+25v+1$ \\
& 3 & $v^6+11v^5+70v^4+148v^3+95v^2+20v+1$ \\
\hline
$E_7$ & 4 & $v^9+13v^8+94v^7+512v^6+1600v^5+2312v^4+1499v^3+395v^2+38v+1$\\
& 5 & $v^8+10v^7+63v^6+292v^5+685v^4+700v^3+260v^2+32v+1$ \\
& 6 & $3v^6+39v^5+172v^4+312v^3+170v^2+27v+1$ \\
\hline
$E_8$ & 10 & $v^9+17v^8+135v^7+719v^6+2568v^5+4652v^4+3014v^3+699v^2+52v+1$\\
& 11 & $v^8+12v^7+92v^6+518v^5+1766v^4+1693v^3+516v^2+46v+1$ \\
& 12 & $5v^6+67v^5+660v^4+964v^3+386v^2+41v+1$
\\ \hline
\end{tabular}
\medskip
\caption{$k(U(q),U^{(l)}(q))$, as polynomials in $v = q-1$} \label{Tab:kuql}
\end{table}

\bigskip

{\bf Acknowledgments}: This research was funded in part by EPSRC
grant EP/D502381/1. The first author thanks Chris Parker for a
discussion about why the groups $U(q)$ are not isomorphic for $G$ of
type $B_r$ and $C_r$.  We would like to thank the referee for many
helpful suggestions.


\bigskip

\begin{thebibliography}{00}

\bibitem{A}
J.~L.~Alperin,
\emph{Unipotent conjugacy in general linear groups},
Comm.\ Algebra \textbf{34}  (2006),  no.\ 3, 889--891.

\bibitem{BH}
H.~B\"urgstein and W.~H.~Hesselink, \emph{Algorithmic orbit
classification for some Borel group actions}, Comp.~Math.
\textbf{61} (1987), 3--41.

\bibitem{C}
R.~W.~Carter, \emph{Finite groups of Lie type. Conjugacy classes and
complex characters}, Pure and Applied Mathematics, New York, 1985.

\bibitem{DM} F.~Digne and J.~Michel, {\em Representations of
finite groups of Lie type}, London Mathematical Society Student
Texts \textbf{21}, Cambridge University Press, Cambridge, 1991.

\bibitem{GAP}
The GAP~Group, \emph{GAP -- Groups, Algorithms, and Programming,
Version 4.3}; 2002, \\ \verb+(http://www.gap-system.org)+.

\bibitem{GDOOBS} S.~M.~Goodwin,
{\em Algorithmic testing for dense orbits of Borel subgroups}, J.\
Pure Appl.\ Algebra \textbf{197} (2005), no.\ 1--3, 171--181.

\bibitem{Gspring}
\bysame,
\emph{Relative Springer isomorphisms}, J.\ Algebra \textbf{290} (2005),
no.\ 1, 266--281.

\bibitem{Gconj} \bysame, 
\emph{On the conjugacy classes in maximal unipotent subgroups of
simple algebraic groups}, Transform.\ Groups \textbf{11} (2006),
no.\ 1, 51--76.

\bibitem{Gzeta} \bysame, 
\emph{Counting conjugacy classes in Sylow $p$-subgroups of Chevalley
groups}, J. Pure Appl.\ Algebra, \textbf{210} (2007), no.\ 1, 201--218.

\bibitem{GRfin} S.~M.~Goodwin and G.~R\"ohrle,
\emph{Finite orbit modules for parabolic subgroups of exceptional
groups}, Indag.\ Math.\ \textbf{15} (2004), no.\ 2, 189--207.

\bibitem{GR} \bysame, 
\emph{Rational points on generalized
flag varieties and unipotent conjugacy in finite groups of Lie
type}, Trans.\ Amer.\ Math.\ Soc., to appear (2008).

\bibitem{H}
G.~Higman, \emph{Enumerating $p$-groups. I. Inequalities}, Proc.\
London Math.\ Soc.\ (3) \textbf{10} (1960), 24--30.

\bibitem{isaacs1}
I.~M.~Isaacs,
\emph{Characters of groups associated with finite
algebras}, J.\ Algebra \textbf{177} (1995), 708--730.

\bibitem{isaacs2}
\bysame,
\emph{Counting characters of upper triangular
groups}, J.\ Algebra \textbf{315} (2007),  no.\ 2, 698--719.

\bibitem{lehrer}
G.~I.~Lehrer,
\emph{Discrete series and the unipotent subgroup},
Compos.\ Math.\ \textbf{28} (1974), 9--19.

\bibitem{malcev}
A.~Mal'cev, \emph{Commutative subalgebras of semi-simple Lie
algebras}, Izv.\ Akad.\ Nauk.\ SSR Ser.\ Math.\ {\bf 9} (1945),
291--300.

\bibitem{robinson}
G.~R.~Robinson, \emph{Counting conjugacy classes of unitriangular
groups associated to finite-dimensional algebras}, J.\ Group Theory
\textbf{1} (1998), no.\ 3, 271--274.

\bibitem{thompson}
J.~Thompson, $k(\U_n(F_q))$, Preprint,
{\tt{http://www.math.ufl.edu/fac/thompson.html}}.

\bibitem{VA}
A.~Vera-L\'opez and J.~M.~Arregi, \emph{Conjugacy classes in
unitriangular matrices}, Linear Algebra Appl.\  \textbf{370} (2003),
85--124.

\end{thebibliography}
\end{document}